# Stationary Transformation of Integrated Brownian Motion

By Eugene Wong
*University of California at Berkeley*

Abstract

Consider an n-fold integrated Brownian motion $W_n(t)$. We show that a simple change in time and scale transforms it into a *stationary* Gaussian process $X_n(t)$. The collection $\{X_n(t)\}$ not only constitutes an interesting family of processes, but their spectral representation is also useful in dealing with integrated Brownian motion. We illustrate this by deriving an explicit representation for the joint density function for $\{W_n(t)\}$ and showing some of its properties.

1. Integrated Brownian Motion

Let $B(t), t \geq 0,$ be a standard Brownian motion and define n-fold *integrated Brownian motion* $W_n(t)$ as follows:

(1.1)
$$W_0(t) = B(t)$$
$$W_n(t) = \int_0^t W_{n-1}(s)\,ds$$

Integration by parts yields the alternative representation

(1.2) $$W_n(t) = \frac{1}{n!}\int_0^t (t-s)^n B(ds)$$

The case of n=1 was apparently first considered by Kolmogorov [3] and is sometimes called the *Kolmogorov diffusion*. The general case was introduced by Shepp [10] and sometimes referred to as *primitives* of Brownian motion.

Integrated Brownian motion has been considered in a number of problem settings (See [4] and [5] for some references). In particular, there is substantial recent interest in integrated Brownian motion in connection with "small ball probabilities." (See for example [1], [6] and [9])

The case of n=1 is of special interest, and explicit results have been found for that case that are not easily generalizable for n>1. For example, McKean [8] solved a winding problem and obtained an explicit expression for the distribution of the time of "first return to zero" for $W_1(t)$. Khoshnevisan and Shi [2] considered the distribution of quadratic functional for $W_1(t)$ and found

(1.3) $$E\left[e^{-\frac{\theta^2}{2}\int_0^1 W_1^2(t)dt}\right] = \left(\frac{2}{\cosh^2\sqrt{\theta/2} + \cos^2\sqrt{\theta/2}}\right)^{1/2}$$

Chen and Li [1] generalized (1.3) for n>1 but their result is substantially more complicated.

In [12] and [13] McKean's result was used to obtain explicit expressions for the distribution of inter-zero intervals of a particular stationary Gaussian process, thereby resolving an issue that had been outstanding for some time ([7], [11]). The particular stationary process was obtained from $W_1(t)$ by a simple scale and time transformation. A similar transformation of $W_n(t)$ yields a family of jointly stationary Gaussian processes that are the principal subject of this study.

2. Some Basic Properties of Integrated Brownian Motion

Because $W_n(t)$ is the (n+1)-fold integral of a white noise, it enjoys a renewal property that can be stated as follows: Define

(2.1) $$\tilde{W}_n(t) = W_n(t+s) - \sum_{k=0}^n \frac{t^k}{k!} W_{n-k}(s)$$

Then $\tilde{W}_n(t)$ is again an n-fold integrated Brownian motion. For the collection $W(t) = (W_0(t),....,W_n(t))$ define its density function

(2.2) $$\pi(w,t) = P\{W(t) \in dw\}/dw \qquad w \in R^{n+1}$$

Then, the renewal property implies that

(2.3) $$\pi_a(w,t) = P(W(t+s) \in dw \mid W(s) = a)/dw = \pi(w - \mu(a,t), t)$$

where

(2.4) $$\mu_n(a,t) = \sum_{k=0}^n \frac{t^k}{k!} a_{n-k}$$

Because $W(t)$ is Markov, all finite dimensional distributions can now be expressed in terms of the static density $\pi(w,t)$.

Since the components of $W(t)$ are jointly Gaussian, we can immediately write, as in [L2],

$$(2.5) \quad p(w,t) = \frac{1}{(\sqrt{2\pi})^n |R(t)|^{1/2}} e^{-\frac{1}{2}w'R^{-1}(t)w}$$

where $w$ is treated as a column vector and $R(t)$ is the covariance matrix with elements

$$(2.6) \quad (R(t))_{jk} = EW_j(t)W_k(t) = \frac{1}{j!k!}\int_0^t (t-s)^{j+k}\,ds$$

$$= \frac{t^{j+k+1}}{j!k!(j+k+1)}$$

We note that the entries of $R(t)$ are independent of its overall dimension. Consider a collection of matrices $\{M_{mn}\}$ indexed by the number of rows and columns. We call such a collection *dimension-free* if

$$(M_{mn})_{jk} = (M_{pq})_{jk} \quad \forall j \leq m \wedge p \text{ and } \forall k \leq n \wedge q$$

The covariance matrices are dimension-free, but their inverses are not. Indeed, finding the inverse $R^{-1}(t)$ in (2.5) one dimension at a time is not only tedious, but the resulting representation is not particularly useful unless one is only interested in a single dimension. It turns out that although $R^{-1}(t)$ is not dimension-free, it can be written as a product of dimension-free matrices. This result will be derived using the spectral representations made possible by transforming integrated Brownian motion into stationary processes.

3. Transforming Integrated Brownian Motion

For $n \geq 0$, define

$$(3.1) \quad X_n(t) = e^{-(n+\frac{1}{2})t} W_n(e^t)$$

We note that

$$(3.2) \quad EX_0(t)X_0(s) = e^{-\frac{1}{2}(t+s)} \min(e^t, e^s) = e^{-\frac{1}{2}|t-s|}$$

Hence, $X_0(t)$ is stationary. Indeed, it is simply the Ornstein-Uhlenbeck process with a particular choice of time scale.

For $n \geq 1,$ we can write

$$X_n(t) = e^{-(n+\frac{1}{2})t}\left(\int_0^{e^t} W_{n-1}(s)ds\right)$$

(3.3)
$$= \int_{-\infty}^{t} e^{-(n+\frac{1}{2})(t-s)} X_{n-1}(s)ds$$

It shows that $X_n$ is a convolution $f_n * X_{n-1}$ where

(3.4) $$f_n(t) = e^{-(n+\frac{1}{2})t} 1(t)$$

where $1(t)$ is the unit step defined by: $1(t) = 1$ for $t \geq 0$ and $= 0$ for $t < 0$.

It follows by induction that $X_n(t)$ is stationary for every n. Indeed, if we denote $X(t) = (X_0(t),....., X_n(t))$ then the components of $X(t)$ are *jointly* stationary and Gaussian. The process $X_n(t)$ might be called an *integrated Ornstein-Uhlenbeck* (IOU) process of order *n*.

Equation (3.3) can be iterated to yield a "white noise" representation

(3.5) $$X_n(t) = \int_{-\infty}^{\infty} h_n(t-s)Z(ds)$$

where Z has independent Gaussian increments with $E[Z(t) - Z(s)]^2 = |t - s|$. We also have a spectral representation

(3.6) $$X_n(t) = \frac{1}{2\pi} \int_{-\infty}^{\infty} H_n(v)e^{ivt} \zeta(dv)$$

where $E\zeta(dv)\bar{\zeta}(dv') = 2\pi\delta_{vv'}dv$. Equations (3.5) and (3.6) have the interpretation that $X_n(t)$ is the result of filtering a Gaussian white noise Z by a linear and time-invariant system with transfer function $H_n$. We denote this fact symbolically by writing

$$X_n = H_n \otimes Z$$

From (3.3) it is easy to see that

(3.7) $$H_n(v) = \prod_{k=0}^{n} \frac{1}{(k+1/2+iv)}$$

One can also invert the Fourier transform to get

(3.8) $$h_n(t) = \frac{1}{n!} e^{-\frac{1}{2}t}(1-e^{-t})^n 1(t)$$

But it is usually easier to work with the transfer function $H_n$ because computation with $H_n$ is mostly algebraic. For example, to compute the cross-correlation functions, we can write

(3.9) $$EX_j(t)X_k(s) = \frac{1}{2\pi} \int_{-\infty}^{\infty} H_j(v)\overline{H}_k(v) e^{iv(t-s)} dv$$

which is easily evaluated using residues.

We summarize the foregoing as follows:

Proposition 3.1

Let $X_n(t)$ be defined by (3.1). Then the processes $X_n(t)$ are jointly Gaussian and stationary. Furthermore, $X_n(t)$ has a spectral representation given by (3.6) with $H_n$ given by (3.7).

4. Spectral Characterization

Let $p_n(x)$ denote the joint density function of $X(t) = (X_0(t),...., X_n(t))$, i.e.,

(4.1) $$p_n(x)dx = P(X(t) \in dx)$$

We can write

(4.2) $$p_n(x_0,...., x_n) = p(x_n | x_0,...., x_{n-1}) p_{n-1}(x_0,...., x_{n-1})$$

Hence, we need only to compute $p(x_n | x_0,...., x_{n-1})$, which is the density for

(4.3) $$\hat{X}_n(t) = X_n(t) - E(X_n(t) | X_0(t),...., X_{n-1}(t))$$

Now, $\hat{X}_n(t)$ is characterized by two properties. First, it is of the form

(4.4) $$\hat{X}_n(t) = X_n(t) + \sum_{k=0}^{n-1} \alpha_{nk} X_k(t)$$

Second,

(4.5) $$E\hat{X}_n(t) X_k(t) = 0 \qquad 0 \le k \le n-1$$

Since $X_k = H_k \otimes Z$ for each $k$, we have $\hat{X}_n = \hat{H}_n \otimes Z$ with

(4.6) $$\hat{H}_n(v) = H_n(v) + \sum_{k=0}^{n-1} \alpha_{nk} H_k(v)$$

It turns out that $\hat{H}_n(v)$ has a simple form.

Proposition 4.1

The process $\hat{X}_n(t)$ has a spectral representation

(4.7) $$\hat{X}_n(t) = \frac{1}{2\pi} \int_{-\infty}^{\infty} \hat{H}_n(v) e^{ivt} \zeta(dv)$$

where $\hat{H}(v)$ is given by

(4.8) $$\hat{H}_n(v) = \frac{n!}{2n!} H_n(v) / \overline{H}_{n-1}(v)$$
$$= \frac{n!}{2n!} \frac{1}{(n+1/2+iv)} \prod_{k=0}^{n-1} \frac{(k+1/2-iv)}{(k+1/2+iv)}$$

*proof:* First, we note that

$$E\hat{X}_n(t) X_k(t) = \int_{-\infty}^{\infty} \hat{H}_n(v) \overline{H}_k(v) dv$$

For $0 \le m \le n-1$, we can write

$$\hat{H}_n(v) \overline{H}_m(v) = \frac{n!}{2n!} \frac{1}{(n+1/2+iv)} \frac{\prod_{k=m+1}^{n-1}(k+1/2-iv)}{\prod_{k=0}^{n-1}(k+1/2+iv)}$$

This can be viewed as the Fourier transform of a convolution $f_1 * f_2$, where

$$f_1(t) = 1(t)e^{-(n+1/2)t}$$

and $f_2$ has a Fourier transform given by

$$F_2(v) = \frac{n!}{2n!} \prod_{k=m+1}^{n-1} \left(\frac{k+1/2-iv}{k+1/2+iv}\right) \prod_{k=0}^{m} \frac{1}{(k+1/2+iv)}$$

Because $F_2$ is an $L_2$ function that is analytic and bounded in the lower-half $v$-plane, $f_2(t)$ is zero for $t < 0$. It follows that for every $m \leq n-1$

$$E\hat{X}_n(t)X_m(t) = \frac{1}{2\pi} \int_{-\infty}^{\infty} \hat{H}_n(v)\overline{H}_m(v)dv$$

$$= \lim_{t \to 0} \int_0^t f_1(t-s)f_2(s)ds = 0$$

The coefficient $\dfrac{n!}{2n!}$ is determined by the $H_n$ term in (4.6). Specifically, we have

$$\left. \hat{H}_n(v) \prod_{k=0}^{n} (k+1/2+iv) \right|_{iv=-(n+\frac{1}{2})} = \frac{n!}{2n!} \prod_{k=0}^{n-1} (k+n+1) = 1$$

as required. ∎

To complete the representation for $\hat{X}_n(t)$, we need to re-express $\hat{H}_n(v)$ as a linear combination of $\{H_k(v), 0 \leq k \leq n\}$ as in (4.6). It is convenient to be free of the coefficient $\dfrac{n!}{2n!}$ and define

(4.9) $$G_n(v) = H_n(v)/\overline{H}_{n-1}(v) = \frac{\prod_{k=0}^{n-1}(k+1/2-iv)}{\prod_{k=0}^{n}(k+1/2+iv)}$$

We now proceed to express $G_n(v)$ as a linear combination of $H_k(v), 0 \le k \le n$.

At this point it is useful to use matrix representation. We adopt the convention that the rows and columns of any matrix are numbered from 0 onward, and define for a matrix $M$ the operation $M^*$ by

(4.10) $\quad (M^*)_{jk} = (-1)^{j+k} (M)_{jk}$

This extends a definition introduced in [5].

We denote by $G(v)$ a column vector with components $G_k(v)$, and similarly $H(v)$. Clearly, $G(v)$ and $H(v)$ are dimension-free in the sense of Section 2 and can be of any size. The relationship between $G(v)$ and $H(v)$ is given as follows:

Proposition 4.2

With $G(v)$ and $H(v)$ defined by (4.9) and (3.7) respectively, we have

(4.11) $\quad G(v) = AH(v)$

where the matrix $A$ and its inverse are given by

(4.12) $\quad \begin{aligned} A_{nm} &= (-1)^{n+m} \frac{(n+m)!}{m!(n-m)!} & n \ge m \\ &= 0 & n < m \end{aligned}$

(4.13) $\quad \begin{aligned} (A^{-1})_{nm} &= \frac{(2m+1)n!}{(n+m+1)!(n-m)!} & n \ge m \\ &= 0 & n < m \end{aligned}$

The derivation of $A$ and $A^{-1}$ is given in the Appendix. We note that both $A$ and $A^{-1}$ are dimension-free in the sense defined in Section 2.

5. Computing Density Functions

The representation of $\hat{X}_n(t)$ in terms of $\{X_k(t)\}$ is given by

$$\hat{H}_n(v) = \left(\frac{n!}{2n!}\right) G_n(v)$$

(5.1)
$$= \left(\frac{n!}{2n!}\right) \sum_{k=0}^{n} (-1)^{n-k} \left(\frac{(n+k)!}{k!(n-k)!}\right) H_k(v)$$

and equivalently

(5.2)
$$\hat{X}_n(v) = \left(\frac{n!}{2n!}\right) \sum_{k=0}^{n} (-1)^{n-k} \left(\frac{(n+k)!}{k!(n-k)!}\right) X_k(v)$$

We note that $E\hat{X}_n^2(t)$ is easily found to be

(5.3)
$$\sigma_n^2 = E\hat{X}_n^2(t) = \frac{1}{2\pi} \int_{-\infty}^{\infty} \left|\hat{H}_n(v)\right|^2 dv = \frac{1}{2n+1}\left(\frac{n!}{2n!}\right)^2$$

Hence, we can write

(5.4)
$$p(x_n \mid x_0,....,x_{n-1}) = \frac{1}{\sqrt{2\pi\sigma_n^2}} e^{-\frac{1}{2\sigma_n^2}(\hat{x}_n)^2}$$

where

$$\hat{x}_n = \left(\frac{n!}{2n!}\right) \sum_{k=0}^{n} A_{nk} x_k$$

and

$$A_{nk} = (-1)^{n+k} \left(\frac{(n+k)!}{k!(n-k)!}\right)$$

From (4.2) we have the following:

Proposition 5.1

The joint density $p_n(x)$ of $(X_0(t),....,X_n(t))$ is given by:

$$p_n(x) = p_0(x_0) \prod_{k=1}^{n} p(x_k \mid x_0, \ldots, x_{k-1})$$

(5.5)
$$= \left( \prod_{k=0}^{n} \frac{1}{\sqrt{2\pi\sigma_k^2}} \right) e^{-\frac{1}{2} \sum_{k=0}^{n} (2k+1) \left( \sum_{j=0}^{k} A_{kj} x_j \right)^2}$$

$$= \left( \prod_{k=0}^{n} \frac{1}{\sqrt{2\pi\sigma_k^2}} \right) e^{-\frac{1}{2} (Ax)' \Lambda (Ax)}$$

where $\Lambda$ is a diagonal matrix with $(\Lambda)_{kk} = (2k+1)$.

We can now evaluate the density of $W(t)$. First, since

$$W_n(t) = t^{n+\frac{1}{2}} X_n(\ln t)$$

we have

$$\pi(w,t) dw = P(W(t) \in dw)$$
$$= P(X(t) \in d\xi(w,t))$$

where

$$\xi_k(w,t) = (t)^{-(k+1/2)} w_k$$

Hence,

(5.6)
$$\pi(w,t) = (t)^{-\frac{1}{2}(n+1)^2} p_n(\xi(w,t))$$

and

(5.7)
$$\pi_a(w,t) = \pi(w - \mu(a,t), t)$$
$$= (t)^{-\frac{1}{2}(n+1)^2} p_n(\xi(w - \mu(a,t), t))$$

We note that

(5.8)
$$\xi_k(w - \mu(a,t),t) = t^{-(k+1/2)}\left[w_k - \sum_{j=0}^{k}\frac{t^j}{j!}a_{k-j}\right]$$
$$= \frac{w_k}{t^{k+1/2}} - \sum_{j=0}^{k}\frac{1}{(k-j)!}\left(\frac{a_j}{t^{j+1/2}}\right)$$

Hence, we can write

(5.9) $\quad \xi(w - \mu(a,t),t) = T(t)w - \Gamma T(t)a$

where $T(t)$ is a diagonal matrix with entries $(T(t))_{kk} = t^{-(k+1/2)}$, and $\Gamma$ is defined by (A.3). Since from the Appendix we have

(5.10) $\quad A\Gamma = A^*$

The transition density can now be written as follows:

Proposition 5.2

The transition density for $W(t)$ can be expressed as

(5.11)
$$\pi_a(w,t) = (t)^{-\frac{1}{2}(n+1)^2} p_n(T(t)w - \Gamma T(t)a))$$
$$= (t)^{-\frac{1}{2}(n+1)^2} K_n e^{-\frac{1}{2}[AT(t)w - A^*T(t)a]'\Lambda[AT(t)w - A^*T(t)a]}$$

where the normalizing constants $K_n$ are given by

$$K_n = \frac{1}{(2\pi)^{(n+1)/2}}\sqrt{\frac{(2n+1)!}{2^n n!}}\left(\prod_{m=0}^{n}\frac{2m!}{m!}\right)$$

It is easy to see from (5.10) that the transition density has the symmetry

$$\pi_a(w,t) = \pi_{w^*}(a^*,t)$$

as shown in [5].

Equation (5.11) shows that $W(t)$ satisfies a modified independent-increments property, viz., for $t_0 < t_1 < .... < t_N$, $\{W(t_n) - B(t_n - t_{n-1})W(t_{n-1}), 1 \le n \le N\}$ are independent vector random variables. The matrix $B$ is given by

$$B(t) = T^{-1}(t)\Gamma T(t)$$

The derivation of the matrix $A$ and its inverse facilitates computation involving integrated Brownian motion in many ways. One example is computing the inverse to the covariance matrix in Section 2. From (5.11) it is clear that the quantity $R^{-1}(t)$ in (2.5) is given by

(5.11) $\qquad R^{-1}(t) = T(t)A'\Lambda A T(t)$

It is useful to compute some numerical examples. From (2.6) we have

(5.12) $\qquad (R^{-1}(t))_{jk} = t^{-(j+k+1)} j!k!(\rho)_{jk}$

where $\rho^{-1}$ is a dimension-free matrix with entries

$$(\rho^{-1})_{jk} = \frac{1}{j+k+1}$$

The matrix $\rho$, however, is neither dimension-free nor simple. For dimension $N$, $\rho$ is given by

$$(\rho)_{jk} = (-1)^{j+k} \sum_{m=j\vee k}^{N} \frac{(j+m)!(k+m)!(2m+1)}{(j!)^2(k!)^2(m-j)!(m-k)!}$$

For $N=1$ and 2, $\rho$ is evaluated to be as follows:

$$N = 1 \qquad \rho^{-1} = \begin{bmatrix} 1 & 1/2 \\ 1/2 & 1/3 \end{bmatrix} \qquad \rho = \begin{bmatrix} 4 & -6 \\ -6 & 12 \end{bmatrix}$$

$$N = 2 \qquad \rho^{-1} = \begin{bmatrix} 1 & 1/2 & 1/3 \\ 1/2 & 1/3 & 1/4 \\ 1/3 & 1/4 & 1/5 \end{bmatrix} \qquad \rho = \begin{bmatrix} 9 & -36 & 30 \\ -36 & 192 & -180 \\ 60 & -180 & 180 \end{bmatrix}$$

and $\rho\rho^{-1} = I$ is easily verified by direct multiplication. We note that the entries of $\rho$ oscillate in sign and grow in absolute value quite rapidly with the overall dimension. It is preferable to work with the representation (5.11) for $R^{-1}(t)$ instead of (5.12).

Appendix

From (4.9) we can write

$$\frac{\prod_{k=0}^{n-1}(k+1/2-iv)}{\prod_{k=0}^{n}(k+1/2+iv)} = \sum_{k=0}^{n} A_{nk} \left( \prod_{j=0}^{k} \frac{1}{(j+1/2+iv)} \right)$$

Multiplying both sides by $\prod_{k=m}^{n}(k+1/2+iv)$ and setting $iv = -(m+1/2)$, we get

$$\frac{(n+m)!}{m!(n-m)!} = \sum_{k=m}^{n} \frac{A_{nk}}{(k-m)!}$$

which can be put in matrix form

(A.1) $\qquad B = A\Gamma$

with

(A.2) $\qquad B_{nm} = \frac{(n+m)!}{m!(n-m)!} \qquad n \geq m$
$\qquad\qquad\quad = 0 \qquad\qquad\qquad n < m$

and

(A.3) $\qquad \Gamma_{nm} = \frac{1}{(n-m)!} \qquad n \geq m$
$\qquad\qquad\quad = 0 \qquad\qquad\qquad n < m$

We now compute

$$(\Gamma\Gamma^*)_{jk} = 0 \qquad j < k$$
$$\qquad\qquad = 1 \qquad j = k$$
$$\qquad\qquad = \sum_{m=k}^{j} \frac{1}{(j-m)!} \frac{(-1)^{m+k}}{(m-k)!} = (1-1)^{j-k} = 0 \qquad j > k$$

Hence,

(A.4) $$\Gamma^* = \Gamma^{-1}$$

and

(A.5) $$A = B\Gamma^*.$$

Now, from (A.2) and (A.3), we have

$$\begin{aligned}
(B\Gamma^*)_{nm} &= \sum_{k=m}^{n} \frac{(-1)^{k+m}(n+k)!}{(k-m)!k!(n-k)!} \\
&= \sum_{j=0}^{n-m} \frac{(-1)^j (n+m+j)!}{j!(n-m-j)!(m+j)!} \\
&= (-1)^m \left( \frac{d^n}{dz^n} \sum_{j=0}^{n-m} \frac{(z)^{n+m+j}}{j!(n-m-j)!} \right)_{z=-1} \\
&= \frac{(-1)^m}{(n-m)!} \left( \frac{d^n}{dz^n} \left[ z^{n+m}(1+z)^{n-m} \right] \right)_{z=-1} \\
&= \frac{(-1)^m}{(n-m)!} \binom{n}{m} \frac{(-1)^n (n+m)!(n-m)!}{n!} \\
&= (-1)^{n+m} \frac{(n+m)!}{m!(n-m)!} = (B^*)_{nm}
\end{aligned}$$

Hence, (A.5) implies that

(A.6) $$A = B\Gamma^* = B^*$$

and (4.12) follows.

To prove (4.13), we note that for each $n$, $H_n(v)$ must be a linear combination of $\{G_k(v), 0 \le k \le n\}$. Hence, we can write

$$H(v) = CG(v)$$

Since

$$\frac{1}{2\pi} \int_{-\infty}^{\infty} G_j(v)\overline{G}_k(v)dv = \frac{\delta_{jk}}{2k+1}$$

we have

$$\frac{(C)_{nk}}{2k+1} = \frac{1}{2\pi}\int_{-\infty}^{\infty} H_n(v)\overline{G}_k(v)dv$$

$$= \frac{1}{2\pi}\int_{-\infty}^{\infty}\left(\prod_{j=k}^{n}\frac{1}{(j+1/2+iv)}\right)\left(\prod_{j=0}^{k-1}\frac{1}{(j+1/2-iv)}\right)dv$$

$$= \sum_{m=k}^{n}\frac{(-1)^{m-k}m!}{(m-k)!(n-m)!(m+k+1)!}$$

Now, let

$$f(z) = \sum_{m=k}^{n}\frac{m!z^{m+k+1}}{(m-k)!(n-m)!(m+k+1)!}$$

Then, we have

$$\frac{(C)_{nk}}{2k+1} = -f(-1)$$

Since

$$\frac{d^{k+1}f(z)}{dz^{k+1}} = \sum_{m=k}^{n}\frac{z^m}{(m-k)!(n-m)!} = \frac{1}{(n-k)!}z^k(1+z)^{n-k}$$

and

$$\left(\frac{d^j f(z)}{dz^j}\right)_{z=0} = 0 \quad 0 \le j \le k$$

we can write

$$f(z) = \frac{1}{k!}\int_0^z (z-y)^k \left(\frac{d^{k+1}f(y)}{dy^{k+1}}\right)dy$$

$$= \frac{1}{k!}\int_0^z (z-y)^k \frac{1}{(n-k)!}y^k(1+y)^{n-k}dy$$

It follows that

$$f(-1) = \frac{(-1)^{k+1}}{k!(n-k)!} \int_{-1}^{0} y^k (1+y)^n dy$$

$$= -\frac{n!}{(n-k)!(n+k+1)!}$$

and

(A.7) $\quad (C)_{nk} = \dfrac{(2k+1)n!}{(n-k)!(n+k+1)!}$

Since $G_k(v)$ are linearly independent functions and

$$G(v) = ACG(v)$$

we must have

(A.8) $\quad A^{-1} = C$

and (4.13) is proved.

Note that we have derived $A$ and $A^{-1}$ separately by means of spectral representation without inverting the matrix $A$.

# References

[1] Chen, X. and Li, W. V. (2003). "Quadratic functional and small ball probabilities for the m-fold integrated Brwonian motion," Ann. Prob. **31** 1052-1077.

[2] Khoshenevisan, D. and Shi, Z. (1998). "Chung's law for integrated Brownian motion," Trans. Am. Math. Soc. **350** 4253-4264.

[3] Kolmogorov, A. (1934). "Zufällige Bewegungen," Ann. of Mat. (2) **35** 116-117.

[4] Lachal, A. (1997). "Local Asymptotic classes for the successive primitives of Brownian motion," Ann. Prob. **25** 1712-1734.

[5] Lachal, A. (1997). "Regular points for the successive primitives of Brownian motion," J. Math. Kyoto Univ. **37** 99-119.

[6] Li, W. V. and Shao, Q. M. (2002). "Gaussian processes: inequalities, small ball probabilities and applications," in Stochastic Processes: theory and Methods. Handbook of Statistics (C.R. Rao and D. Shanbhag, eds.) **19** 533-598. North Holland, Amsterdam.

[7] Longuet-Higgins, M. S. (1963). "Bounding approximations to the distribution of intervals between zeros of a stationary Gaussian process, " in Time Series Analysis (M. Rosenblatt ed.) 63-88. John Wiley, New York.

[8] McKean, H. P., Jr. (1963). "A winding problem for a resonator driven by a white noise," J. Math. Kyoto Univ. **2** 227-235.

[9] Nazarov, A. I. and Nikitin, Y. Y. (2004). "Exact $L_2$-small ball behavior of integrated Gaussian processes and spectral assymptotics of boundary value problems," Prob. Th. and Related Fields **129** 469-494.

[10] Shepp, L. A. (1966). "Radon-Nikodym derivatives of Gaussian measures," Ann. Math. Stat. **37** 321-354.

[11]     Slepian, D. (1963). "On the zeros of Gaussian noise," in Time Series Analysis, op. cit., 104-115.

[12]     Wong, E. (1966). "Some results concerning the zero-crossings of Gaussian noise,' SIAM J. Appl. Math. **14** 1246-1254.

[13]     Wong, E. (1970). "The distribution of intervals between zeros for a stationary Gaussian process," SIAM J. Appl. Math. **18** 67-73.